\input{aipcheck}

\documentclass[
    ,final            
  ]
  {aipproc}

\layoutstyle{8x11single}

\usepackage[english]{babel}
\usepackage{amsmath}
\usepackage{amssymb}

\begin{document}

\title{An Example of Symmetry Exploitation for Energy-related Eigencomputations}


\author{M.~Petschow}{
  address={RWTH Aachen, AICES, Aachen, Germany}
}

\author{E.~Di Napoli}{
  address={RWTH Aachen, AICES, Aachen, Germany}
}

\author{P.~Bientinesi}{
  address={RWTH Aachen, AICES, Aachen, Germany}
}


\begin{abstract}
One of the most used approaches in simulating materials is the tight-binding approximation. When using this method in a material simulation, it is necessary to compute the eigenvalues and eigenvectors of the Hamiltonian describing the system. In general, the system possesses few explicit symmetries. Due to them, the problem has many degenerate eigenvalues. The ambiguity in choosing a orthonormal basis of the invariant subspaces, associated with degenerate eigenvalues, will result in eigenvectors which are not invariant under the action of the symmetry operators in matrix form. A meaningful computation of the eigenvectors needs to take those symmetries into account. A natural choice is a set of eigenvectors, which simultaneously diagonalizes the Hamiltonian and the symmetry matrices. This is possible because all the matrices commute with each other. The simultaneous eigenvectors and the corresponding eigenvalues will be in a parametrized form in terms of the lattice momentum components. This functional dependence of the eigenvalues is the dispersion relation and describes the band structure of a material. Therefore it is important to find this functional dependence in any numerical computation related to material properties.
\end{abstract}

\maketitle


\section{Introduction}

Tight-binding (TB) is a method used to investigate the
electronic structure of a large class of solid materials~\cite{sla54}. When used
in conjunction with numerical simulations{\bf ,} this method introduces several
simplifications that reduce the complexity of the description of the material.
Every solid material is constituted of atomic nuclei that identify a
lattice and are the source of potential energy. On other hand,
the nucleus-nucleus and electron-electron interactions are neglected. The
TB model assumes that the electrons are tightly bound to their corresponding
nuclei, implying that their wave functions are localized. Furthermore{\bf ,} atoms
interact weakly only through their valence electrons.

Since the electrons are moving independently{\bf ,} the Hamiltonian $H$ of the system
is given as a sum of the kinetic energies of the electrons $(p_i^2/2m)$ and the potentials
due to the nuclei $\sum_i V(r_i - R_n)$, with $R_n$ being the positions of the nuclei in three-dimensional space $r_i$. Thus:
\begin{equation}
 H = \sum_{i=1}^{N_e} H_i = \sum_{i=1}^{N_e} \left( \frac{p_i^2}{2m} + \sum_{n=1}^{N_s} V(r_i - R_n) \right) \,, \label{hamilton}
\end{equation}
where $N_e$ denotes the number of considered electrons and $N_s$
the number of lattice sites of the crystal~\cite{czy08,gro00}.
To find the eigenstates of this system{\bf ,} a linear combination of the atomic
orbitals (LCAO) $\sum_i \tilde{v}_i \phi_i$ is used as an ansatz. The atomic orbitals $\phi_i$ are the eigenstates
of the Hamiltonian for an isolated
atom and $\tilde{v}_i$ the coefficients to be computed with the constraint that $|\sum_i \tilde{v}_i \phi_i|=1$.
Since the overlap of the atomic orbitals of neighboring atoms is assumed to be
small{\bf ,} they are treated as orthonormal, i.e. their inner products are $(\phi_i, \phi_j) = \delta_{ij}$. Using this property and the LCAO as an
ansatz{\bf ,} one obtains the following eigenproblem:
\begin{equation}
 H v_n = e_n v_n \quad\mbox{with}\quad n=0,1,\ldots,K-1 \,, \label{eigproblem}
\end{equation} 
where $H \in \mathbb{C}^{K \times K}$ is the Hamiltonian in the basis of the
atomic orbitals. The quantity $v_n \in \mathbb{C}^K$ is a vector of coefficients $\tilde{v}_i$ of the LCAO and the eigenvalue $e_n \in \mathbb{R}$ is the associated energy level. Note that the Hamiltonian
is hermitian and therefore the eigenvalues are real.

The entries of the Hamiltonian are given by
\begin{equation}
 H_{k \ell} = \delta_{k\ell} \alpha_k + \beta_{k\ell} \,,
\label{Hentries}
\end{equation}
where $\beta_{k\ell}$ is the result of an overlap integral between neighboring
electronic orbitals and the underlying lattice potential~\cite{gro00}. In the
simplest case of equal atoms and only nearest neighbor interaction, the expressions in Eq.~(\ref{Hentries}) simplify to $\beta_{k\ell} = -t \left(
  \delta_{k,\ell+1} + \delta_{k+1,\ell} \right)$ and $\alpha_k = \alpha$, where
$\alpha$ and $t$ are constants~\cite{gro00}. Since $t$ represents the interaction between neighboring atoms, it is often called the hopping term.\\

The quantum mechanical problem of finding the electron wave function is
therefore reduced to the solution of a finite dimensional eigenproblem. Having
computed the eigenvalues and eigenvectors, we aim at expressing them in terms of the
lattice momentum components $\textbf{k} = \left( k_1, k_2,
k_3\right)$. Eventually{\bf ,} the whole set of eigenvalues can be seen as a function
of $\textbf{k}$, called the dispersion relation. This is an important relation from which
we can determine a large set of physical properties of a material \cite{czy08,gro00}. Therefore determining this relation
numerically is our final goal.

\section{A 2-dimensional Example}

\begin{center}
   \begin{figure}[htb]
    \centering
    \includegraphics[scale=.39]{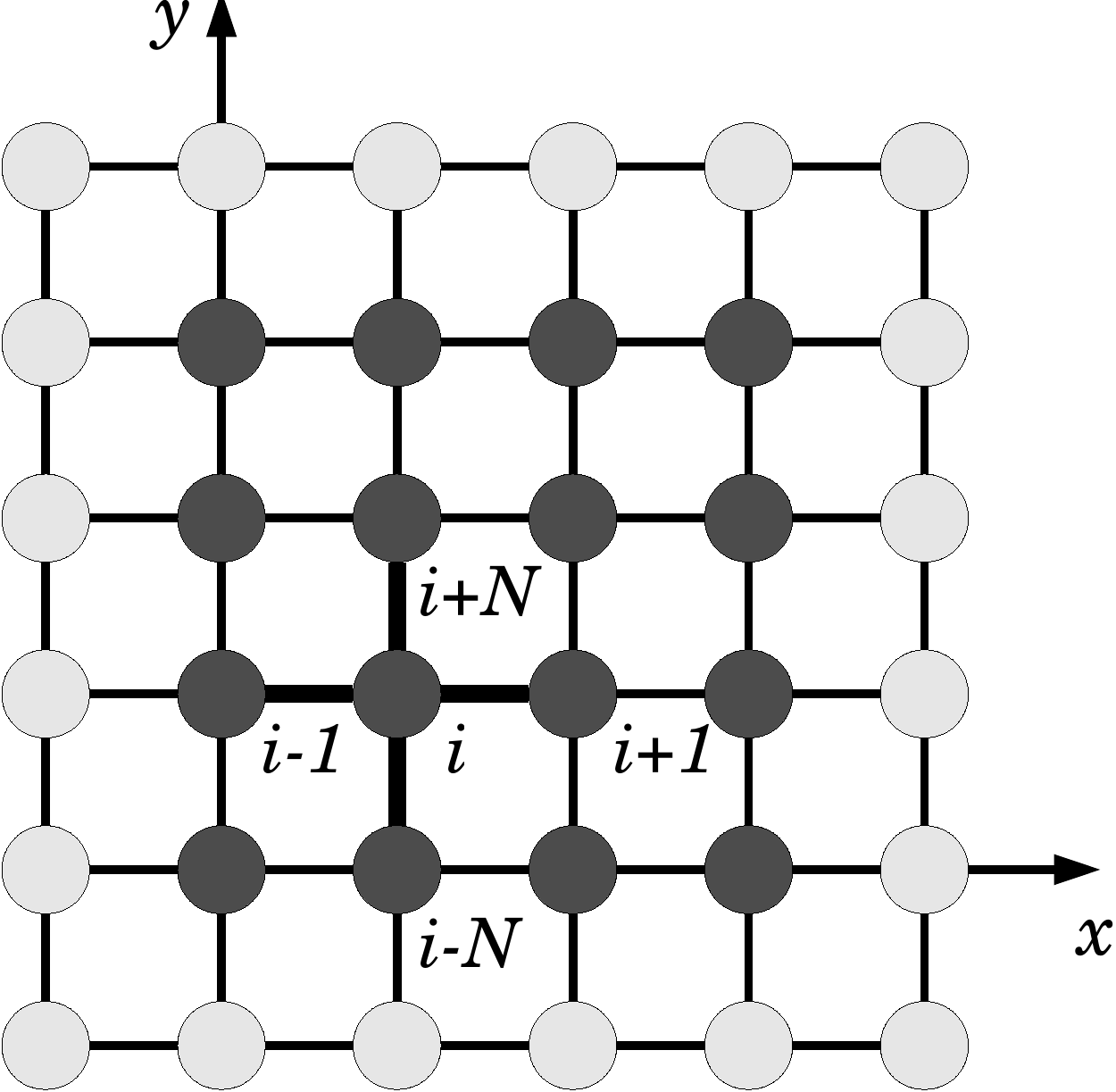} \quad\quad \includegraphics[scale=.3,trim= 0mm 60mm 0mm 70mm]{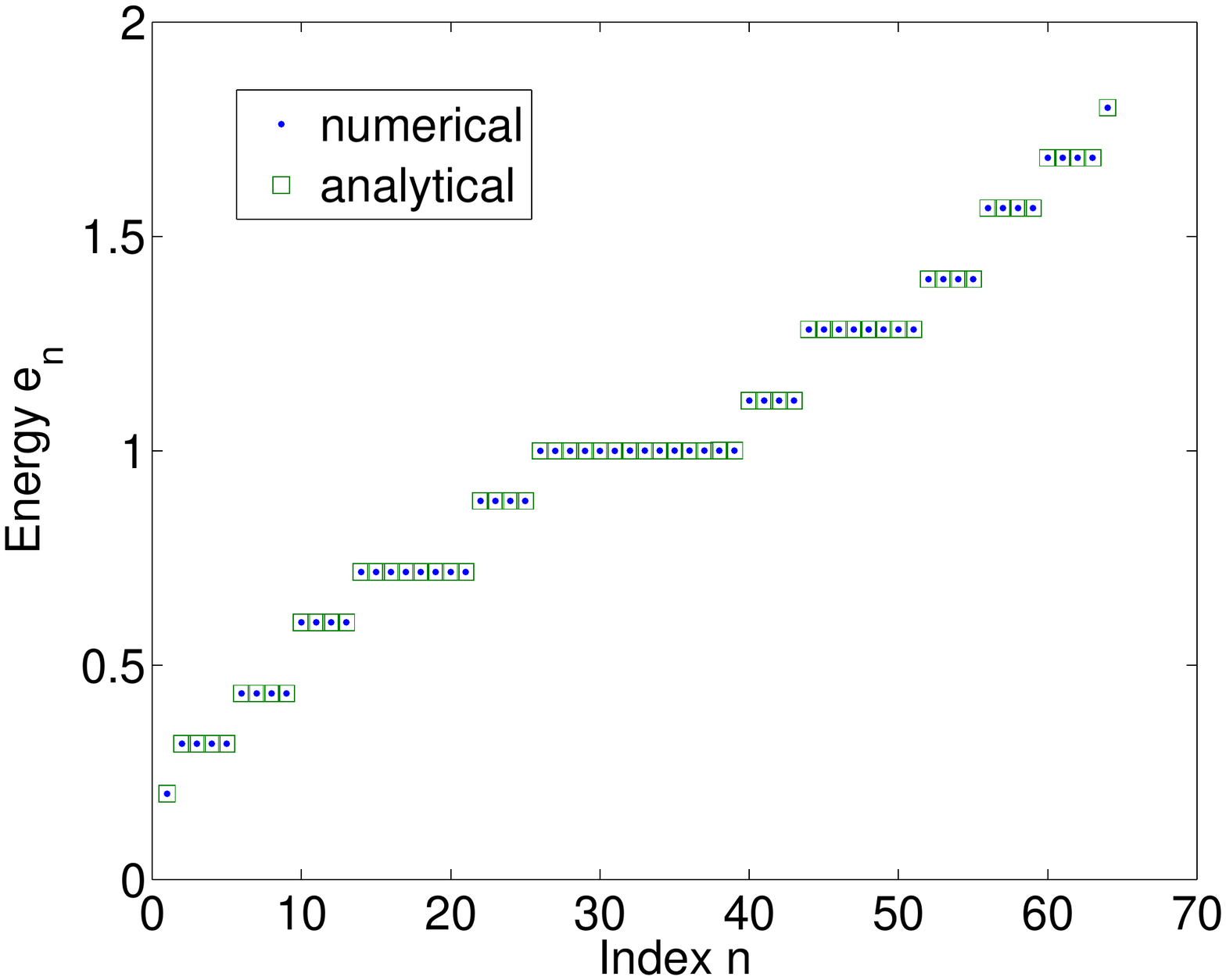}
    \caption{\textit{Left:} Two-dimensional rectangular lattice of equal atoms and nearest neighbor interaction. In this figure the mesh has the size $N=4$; \textit{right:} Spectrum of the Hamiltonian with $\alpha=1.0$ and $t=0.2$ for $N=8$.}
    \label{mesh2d}
   \end{figure} 
\end{center}

In this section we construct a simple example. While it can be solved analytically, we show that it can also be accurately solved numerically. Consider a two-dimensional rectangular lattice of equal atoms, as shown in Figure \ref{mesh2d} (left). The dark colored atoms constitute our $N$-by-$N$ lattice structure and the brighter atoms represent the use of periodic boundary conditions. Each atom in the structure interacts with its four nearest neighbors. The interaction is given by the hopping term $t$ as discussed above.

The Hamiltonian $H \in \mathbb{R}^{N^2 \times N^2}$ of the system has the form
\begin{equation}
 H = \left[ 
       \begin{array}{ccccc}
          C & D &        &        & D \\
          D & C & D      &        &   \\
            & D & C      & \ddots &   \\
            &   & \ddots & \ddots & D \\
          D &   &        & D      & C
       \end{array}
    \right] \,.
    \label{Hform}
\end{equation} 
The matrix $C \in \mathbb{R}^{N \times N}$ is circulant. It is equivalent to the Hamiltonian for the one-dimensional lattice of $N$ identical atoms with periodic boundary conditions and nearest neighbor interactions. It has the same structure as $H$ in Eq.~(\ref{Hform}) with $C$ and $D$ replaced by the scalars $\alpha$ and $-t$, respectively. The matrix $D \in \mathbb{R}^{N \times N}$ is diagonal with all elements equal to $-t$.

The eigenvalues and eigenvectors of $H$ can be expressed in closed form (see \cite{dav79, tee07}). Since $H$ is block-circulant with circulant symmetric blocks, its normalized eigenvectors $v_n \in \mathbb{C}^{N^2}$ are
\begin{equation}
 v_n(r,s) = \frac{1}{N} \left[ 
       \begin{array}{c}
           w_s \\
           \rho_r w_s  \\
           \rho_r^2 w_s  \\
           \vdots  \\
           \rho_r^{N-1} w_s
       \end{array}
    \right] , \quad \mbox{with} \quad
 w_s = \left[ 
       \begin{array}{c}
           1 \\
           \xi_s  \\
           \xi_s^2  \\
           \vdots  \\
           \xi_s^{N-1}
       \end{array}
    \right] \,, \quad
    \rho_r = \exp\left(i\frac{2\pi}{N} r\right) \,, \quad
    \xi_s  = \exp\left(i\frac{2\pi}{N} s\right) \,,
\label{eigvectors}
\end{equation}
where $r,s=0,1,\ldots,N-1$. The parameters $\rho_r$ and $\xi_s$ are the $N$-th roots of unity. The eigenvectors $v_n$ form an orthonormal basis \cite{dav79}. The corresponding eigenvalues are
\begin{equation}
   e_n(r,s) = \alpha - 2t\cos\left(\frac{2\pi}{N}r\right) - 2t\cos\left(\frac{2\pi}{N}s\right) \,.
   \label{eigvals}
\end{equation}

All eigenpairs $(e_n, v_n)$ are parametrized by the quantities $r$ and $s$. The index $n$ of the pair can be defined as any bijective function $n = f(r,s)$.

When defining $k_x$ and $k_y$ as $k_x := 2\pi r/N$ and $k_y := 2\pi s/N$, respectively, Eq.~(\ref{eigvals}) describes the dispersion relation $e(k_x,k_y)$. This relation yields all the allowed energies for possible momenta $k = \sqrt{k_x^2+k_y^2}$ . Because of the periodicity of the crystal, both energy and momentum are quantized \cite{czy08,gro00}. As $N\rightarrow \infty$, the dispersion relation reveals the band structure of the crystal. Therefore it is important to identify this relation. In more complicated cases for which no analytical solution is available, it is important to compute the dispersion relation through a numerical procedure.

In general a numerical computation does not result in eigenvectors in the form of Eq.~(\ref{eigvectors}). For every $m$-degenerate eigenvalue, the matrix $H$ only defines an $m$-dimensional invariant subspace. For instance, from Eq.~(\ref{eigvals}) it can be seen that for even $N$ only the largest and smallest eigenvalues are distinct. All the other eigenvalues have at least multiplicity four. Figure \ref{mesh2d} (right) illustrates this behavior showing the spectrum of a Hamiltonian for $N=8$. The degeneracy creates an ambiguity in choosing a basis for the associated invariant subspaces.

The best a general algorithm can do given only the matrix $H$, is to compute an arbitrary orthonormal basis. Such a solution would not have a parametric expression in terms of the momentum because it does not respect the symmetry of the lattice. In physical terms, it means that the set of eigenvectors is not invariant under the action of the symmetries expressed in matrix form. An obvious solution is to find a vector basis that simultaneously diagonalizes the Hamiltonian $H$ and the symmetries $S_i$. This is indeed possible, because $H$ commutes with the symmetries\footnote{This is a well-known result from Hamiltonian dynamics: the time-dependence of a generic operator $A$ is described by the Eq. $\frac{dA}{dt} \propto [H,A]$. If $A$ represents a conserved symmetry, its derivative with respect to time is automatically null from which the thesis follow.}~\cite{czy08}, generating a closed algebra under multiplication\footnote{The Jacobi identity is readily verified.}:
\begin{equation}
 \left[ H, S_i \right] = 0 \,, \quad \left[ S_i, S_j \right] = 0 \,.
\end{equation}

Our example has two translation symmetries, one along the $x$-axis $S_x = I_N \otimes C_p$, and another along the $y$-axis $S_y = C_p \otimes I_N$. The symbol $\otimes$ denotes the Kronecker or tensor product, $I_N \in \mathbb{R}^{N \times N}$ is the identity matrix, and $C_p \in \mathbb{R}^{N \times N}$ is a circulant matrix with $C_p = \mbox{circ}(0,\ldots,0,1)$. Since both matrices are simply permutation matrices, their inverses are  $S_i^{-1}=S_i^T$, which correspond to translations along the negative direction of the symmetry axis.

Unfortunately, neither $S_x$ nor $S_y$ have distinct eigenvalues, in fact every eigenvalue has multiplicity $N$, and we can only determine invariant subspaces. In order to find the simultaneous eigenvectors of all the matrices, we can create linear combinations that are part of the algebra. Using enough combinations will allow us to identify uniquely (up to a phase) the simultaneous eigenvectors in their parametric form.

For example, we can compute the eigenvectors of $H(S_x-S_y)$ and $S_x(H-S_y)$, and select from their eigenvectors only those that are simultaneous eigenvectors of all the matrices. In this way, we generate a set of $N^2$ simultaneous eigenvectors. They coincide with the analytical solutions of Eq.~(\ref{eigvectors}) up to a phase factor. In our example we can normalize the eigenvectors imposing the first element to be real. The resulting basis is uniquely defined and satisfies the symmetries of the problem.

To show a concrete numerical example, we now look at the results for a 25-by-25 lattice with $\alpha=1.0$ and $t=0.2$: using $\mbox{MATLAB}^\circledR$ for the computation, and denoting computed quantities with a hat, the maximum error of the computed eigenvector entries compared to the analytical solution for real and imaginary part is about $3.5\cdot10^{-14}$; the maximum residual $\mbox{max}_i \parallel H\hat{v}_i - \hat{e}_i \hat{v}_i \parallel_2 \approx 1.9\cdot10^{-13}$; the orthogonality of the computed eigenvectors is $\mbox{max}_{i,j} |\hat{v}_i^*\hat{v}_j - \delta_{ij}| \approx 3.8\cdot10^{-13}$; finally the maximum error in the eigenvalues is $\max_i |\Re(\hat{v}_i^*H\hat{v}_i) - e_i| \approx 1.1\cdot10^{-15}$.

\begin{center}
   \begin{figure}[htb]
    \centering
    \includegraphics[scale=.33,trim= 0mm 60mm 0mm 70mm]{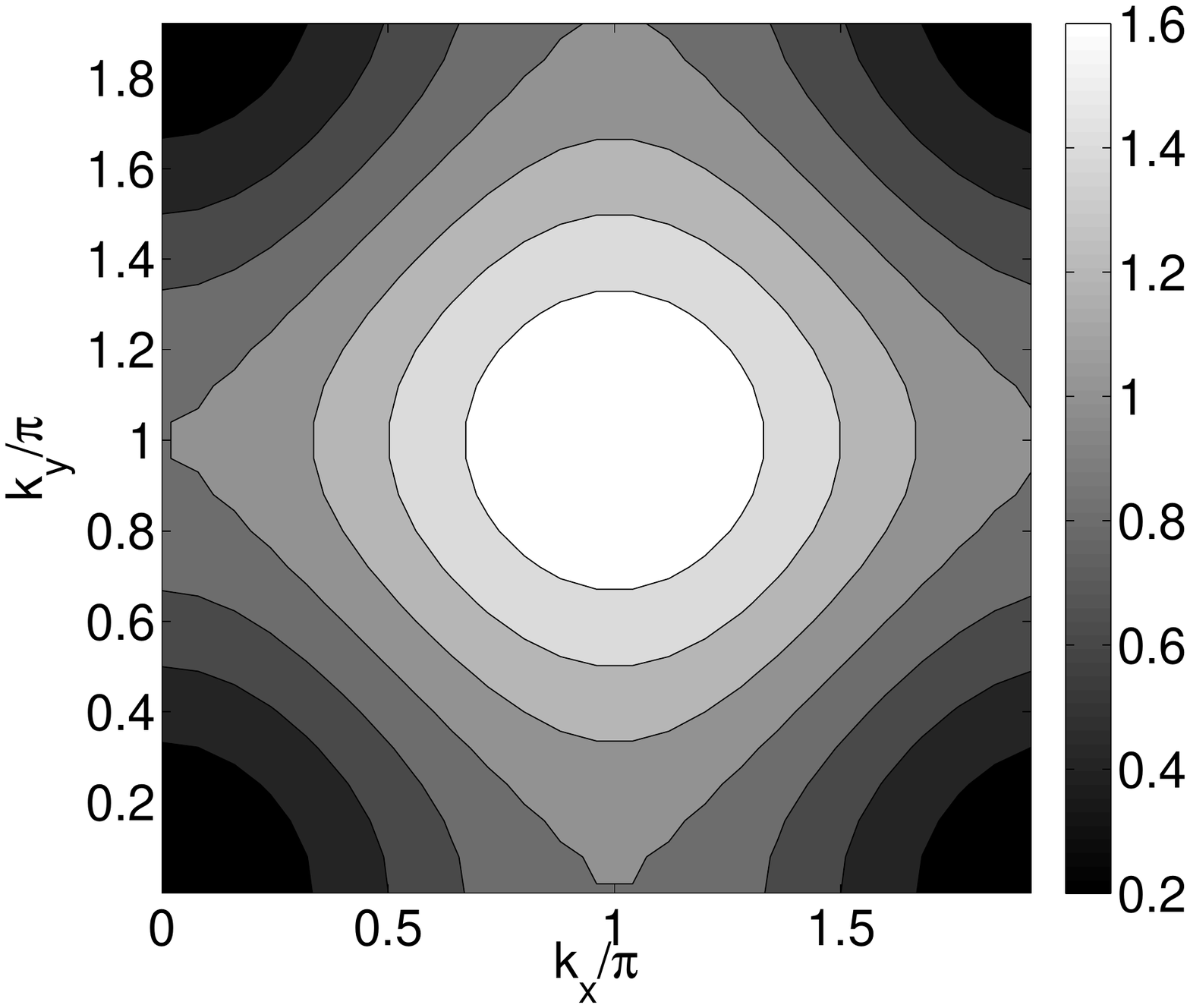} \quad\quad \includegraphics[scale=.33,trim= 0mm 60mm 0mm 70mm]{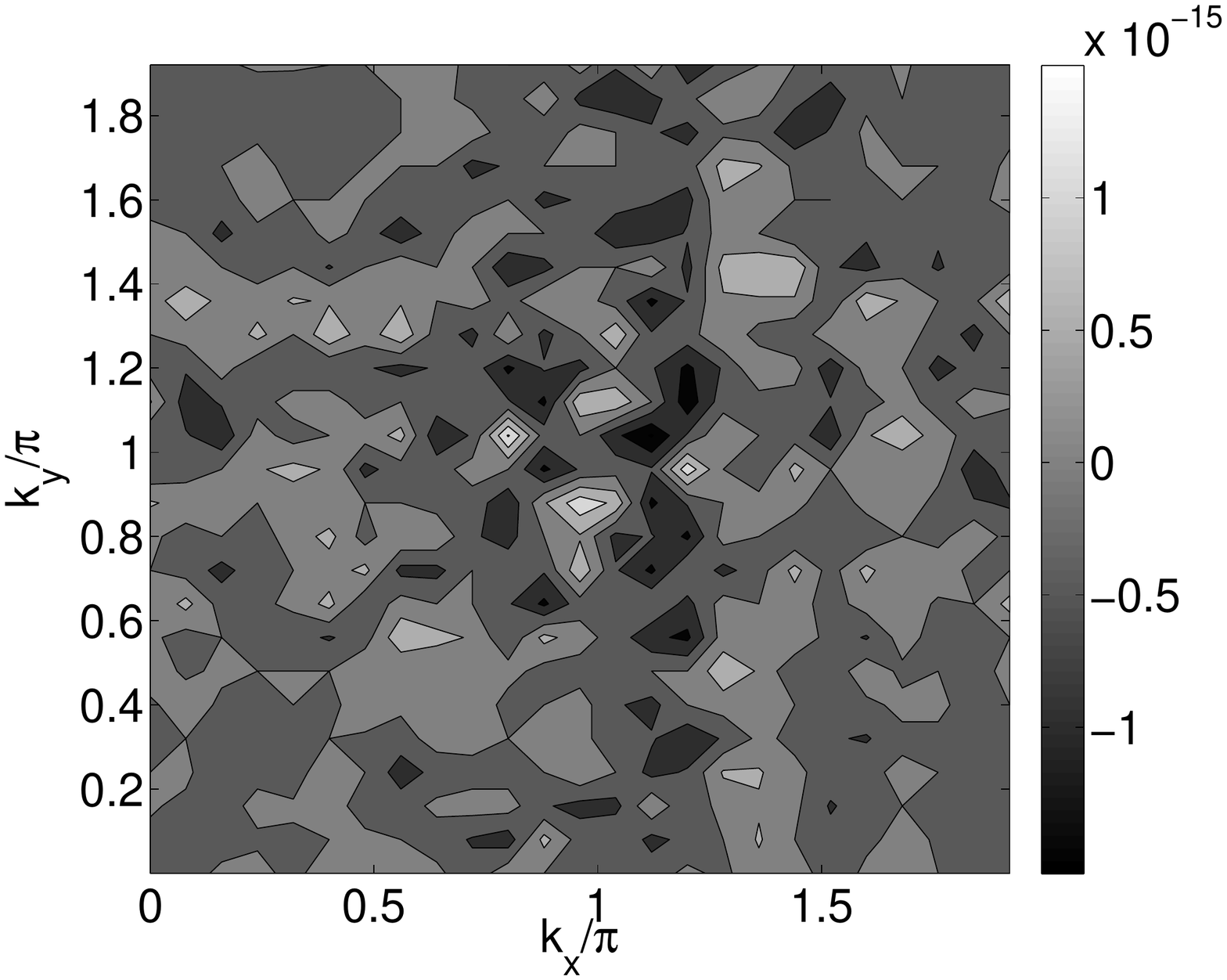}
    \caption{\textit{Left:} The dispersion relation $e(k_x,k_y)$ computed numerically; \textit{right:} The error compared to the analytical solution given by Eq.~(\ref{eigvals}).}
    \label{disrel}
   \end{figure} 
\end{center}

Figure \ref{disrel} (left) shows the dispersion relation computed numerically, while on the right we present the error compared to the analytical solution given by Eq.~(\ref{eigvals}).

\section{Summary and Conclusion}
The computation of eigenvalues and eigenvectors of an Hamiltonian describing a quantum mechanical system can lead to eigenvectors that are not satisfying certain physical requirements.

Through a simple model of a solid material, we show that symmetries introduce degeneracies in the spectrum of the Hamiltonian. These degeneracies lead to an ambiguity in computing a basis for the invariant subspaces associated with the corresponding eigenvalues. A basis, if not chosen carefully, in general will lead to eigenstates which are not taking the symmetries of the problem into account. In order to generate a satisfactory basis the eigenvectors must simultaneously diagonalize the Hamiltonian and the symmetry operators. The eigenvectors that fulfill such a condition compose a complete orthonormal eigenvector basis that is uniquely defined. Finding this basis is the first step in computing the dispersion relation of the material under investigation.

We explored the concept of computing numerically the dispersion relation in a simple model having analytical solutions. Our final goal is to apply this approach to the investigation of irregular materials, where analytical solutions are not known.


\vspace{-3mm}
\begin{theacknowledgments}
Financial support from the Deutsche Forschungsgemeinschaft (German Research Association) through gran GSC 111 is gratefully acknowledged.
\end{theacknowledgments}



\bibliographystyle{aipproc}   


\IfFileExists{\jobname.bbl}{}
 {\typeout{}
  \typeout{******************************************}
  \typeout{** Please run "bibtex \jobname" to optain}
  \typeout{** the bibliography and then re-run LaTeX}
  \typeout{** twice to fix the references!}
  \typeout{******************************************}
  \typeout{}
 }


\end{document}